\documentclass{conm-p-l}
\usepackage{amssymb}
\usepackage[matrix,arrow]{xy}
\numberwithin{equation}{section}
\newtheorem{thm}{Theorem}[section]
\newtheorem{lem}[thm]{Lemma}
\newtheorem{prop}[thm]{Proposition}

\newtheorem{cor}[thm]{Corollary}

\newcommand\End{\mbox{\rm End\,}}

\newcommand\cat{\mbox{\rm cat\,}}

\newcommand\cover[1]{\widetilde{#1}}
\newcommand\complete[1]{\widehat{#1}}
\newcommand\pcomplete[1]{\widehat{#1}_{p}}
\newcommand\apcomplete[1]{\widehat{#1}^{a}_{p}}

\newcommand\field{{\mathbb F}}
\newcommand\integer{{\mathbb Z}}

\newcommand{\incl}{\operatorname{in}}
\newcommand{\homotopic}{\operatorname{\simeq}}
\newcommand\comp{\smash{\lower-.2ex\hbox{\scriptsize$\circ$}}}
\title{A $p$-complete version of the Ganea Conjecture for co-H-spaces}
\author{J. R. Hubbuck and Norio Iwase}
\address[J.~R.~Hubbuck]{Department of Mathematical Sciences, University of Aberdeen, Aberdeen AB24 3QY, United Kingdom}
\email{j.hubbuck@maths.abdn.ac.uk}
\address[N.~Iwase]{Faculty of Mathematics, Kyushu University, Fukuoka 810-8560, Japan}
\email{iwase@math.kyushu-u.ac.jp}
\thanks{This research was partially supported by Grant-in-Aid for Scientific Research (C)08640125 from The Ministry of Science, Sports and Culture.}
\keywords{co-H-space, Ganea's conjecture, almost completion, atomic space}
\subjclass{Primary 55P45}
\date{\today}
\begin{document}
%\maketitle
%
%     ABSTRACT
%
\begin{abstract}
A finite connected CW complex which is a co-H-space is shown to have the homotopy type of a wedge of a bunch of circles and a simply-connected finite complex after almost $p$-completion at a prime $p$.
\end{abstract}
%
%       MANUSCRIPT DOCUMENT
%
\maketitle
\section{Fundamentals and Results}
\par\noindent
When $Y$ is a homotopy associative H-space or when $X$ is a (homotopy) associative co-H-space, the set of based homotopy classes $[Z,Y]$ or $[X,Z]$ respectively, is a group natural in the $Z$ argument.
If the H-multiplication on $Y$ is not known to be homotopy associative, the induced structure on $[Z,Y]$ is that of an algebraic loop; in particular, left and right inverses exist but they may be distinct.
One cannot make as general a statement for $[X,Z]$ when $X$ is a co-H-space.
The immediate problem is that whereas the shearing map for an H-space induces isomorphisms of {\it homotopy} groups, the co-shearing map for a co-H-space induces isomorphisms of {\it homology} groups.
This general situation has been well understood for some decades.

We assume that spaces have the homotopy types of CW-complexes, are based and that maps and homotopies preserve base points.
A space $X$ is a co-H-space if there is a comultiplication map $\nu : X \to X \vee X$ satisfying $j{\comp}\nu \homotopic \Delta : X \to X{\times}X$ where $j : X{\vee}X \hookrightarrow X{\times}X$ is the inclusion and $\Delta$ the diagonal map.
Equivalently, $X$ is a co-H-space if the Lusternik-Schnirelmann category $\cat{X}$ is one.

Statements $(1.1)$ and $(1.2)$ below were shown to be equivalent in \cite{HMR:co-loop}, see also Theorem 0.1 of \cite{Iwase:counter-ganea}.
\begin{enumerate}
\item[(1.1)]\label{(1.1)}~
$X$ is a co-H-space and the comultiplication can be chosen so that $[X,Z]$ is an algebraic loop for all $Z$.
\item[(1.2)]\label{(1.2)}~
The space $X$ has the homotopy type of a wedge of a bunch of circles and a simply-connected co-H-space.
\end{enumerate}

Problem 10 in \cite{Ganea:conjecture} asks ``Is any (non-simply-connected) co-H-space of the homotopy type of $S^1 \vee \cdots \vee S^1 \vee Y$ where there may be infinitely many circles and $\pi_1(Y)=0$?''
The positive statement has become known as `the Ganea conjecture' for co-H-spaces (see section 6 of \cite{Arkowitz:co-H-space}).
The conjecture was resolved thirty years after being raised when the second author constructed in \cite{Iwase:counter-ganea} infinitely many finite complexes which are co-H-spaces but which do not have the homotopy type described in Problem 10.
This leaves open a $p$-complete version of Ganea's conjecture, and probably more difficult, a $p$-local version (see Conjecture 1.6 of \cite{Iwase:counter-ganea}).
The rational version was established in \cite{Henn:almost-rat} as a prime decomposition theorem for an `almost rational' co-H-space.
In this note, we address the $p$-complete problem at a prime $p$.

Some comments are required on $p$-completion.
The $p$-completion of a simply connected co-H-space is rarely a
co-H-space (unless it is a `finite torsion space' in the sense of \cite{CHX:p-atomic}) as a wedge of $p$-complete spaces need not be $p$-complete; it becomes a co-H-object in a categorical sense, which is adequate for some purposes.
More seriously, we are interested in non-simply-connected co-H-spaces and it is shown in \cite{EG:cat1} that when $X$ is a co-H-space, 
\begin{enumerate}
\setcounter{enumi}{2}
\item[(1.3)]\label{(1.3)}~
$\pi_1(X)$ is a free group.
\end{enumerate}
We therefore use fibrewise $p$-completion which we describe after introducing more notation.

Let $X$ be a co-H-space and $B=B\pi_1(X)$, so that $B$ can be chosen as a bunch of circles by $(1.3)$.
Let $i : B \to X$ represent the generators of $\pi_1(X)$ associated with the circles and $j : X \to B$ be the classifying map of the universal cover $p : \cover{X} \to X$.
We may assume that $j{\comp}i$ $\homotopic$ $1_B$, and so $B$ is a homotopy retract of $X$.
Also let $c : X \to C$ be the cofibre of $i : B \to X$, so $C$ is simply connected.
One seems tantalizingly close to Ganea's original conjecture as there are homology equivalences $X \to B \vee C$ and $B \vee C \to X$ inducing isomorphisms of fundamental groups.

For each prime $p$, we consider the fibrewise $p$-completion of $j : X \to B$, $\apcomplete{e} : X \to \apcomplete{X}$ which commutes with projections to $B$.
The map $\apcomplete{e}$ induces an isomorphism of fundamental groups and acts as standard $p$-completion on the fibre $\cover{X}$ and so $\cover{\apcomplete{X}} \homotopic \pcomplete{\cover{X}}$.
%(To simplify notations, we omit the standard circumflex to denote $p$-completion as we do not use $p$-localization nor integral completion.   Thus the standard $p$-completion of $Z$ is denoted by $\pcomplete{Z}$).
Following earlier authors, we refer to this fibrewise $p$-completion for $j : X \to B$ as `almost $p$-completion'.
A general reference for fibrewise $p$-completion is
\cite{BK:completion-localization}, \cite{Bendersky:partial-localization} or \cite{May:fibrewise-local}.
%but a discussion of $p$-localization in our context is given in section 3 of \cite{Iwase:counter-ganea} and this can be adopted to $p$-completion.
Also it is shown in \cite{Iwase:counter-ganea} that a co-H-space $X$ is a co-H-space over $B$ up to homotopy, and so $\apcomplete{X}$ becomes a co-H-object over $B$ in the sense of \cite{James:over-B} and \cite{CHX:p-atomic}.

The main result of this note is the following.

\begin{thm}\label{thm:main}
Let $X$ be a finite, connected, based CW-complex and a co-H-space.
After almost $p$-completion, $\apcomplete{X}$ has the homotopy type of $\apcomplete{(B{\vee}C)}$ where $B$ is a finite bunch of circles and $C$ is a simply-connected finite complex and a co-H-space.
\end{thm}

Let $Y \to B$ be a fibration with cross-section.

\begin{cor}\label{cor:main}
The homotopy set $[\apcomplete{X},\apcomplete{Y}]_B$ inherits an algebraic loop structure from $C$.
\end{cor}

Since $C$ is a simply-connected co-H-space, results of \cite{Xu:thesis} imply that $\pcomplete{C}$ can be decomposed, uniquely up to homotopy, as a completed wedge sum of simply-connected $p$-atomic spaces.

\begin{cor}\label{cor:main2}
$\apcomplete{X}$ has the homotopy type of the almost $p$-completion of a wedge sum of circles and simply-connected $p$-atomic spaces.
%$\apcomplete{X} \homotopic \apcomplete{(\bigvee^{h}_{i=1} S^1 \vee {\bigvee}^{k}_{j=1}{C_{j}})}$, where ${C_j}$ is a simply-connected $p$-atomic space.
\end{cor}

The general strategy used to prove Theorem \ref{thm:main} first occurs in \cite{ISS:hom_co-h} in establishing the Ganea conjecture for complexes of dimension less than $4$.
The existence of a co-H-multiplication enabled the authors to construct a new splitting $C \to X$ to obtain a homotopy equivalence $B \vee C \to X$.
In our case, we adopt techniques for simply-connected $p$-complete spaces of \cite{Hubbuck:p-atomic} and \cite{Xu:thesis} for a similar purpose.

The first author expresses his gratitude to the Kyushu University and the second to the University of Aberdeen for their hospitality.
\section{Finite co-H-complexes}\label{sect:proof1}
\par
We give a different proof of Theorem 3.1 of \cite{ISS:hom_co-h} and include a converse statement for completeness.

\begin{thm}\label{thm:dominate}
Let $X$ have the homotopy type of a based CW complex.
Then $X$ has the homotopy type of a finite complex which is a co-H-space if and only if there are a finite bunch of circles, a connected finite complex $D$ and maps
$$
\rho : B{\vee}{\Sigma{D}} \to X\quad\text{and}\quad\sigma : X \to B{\vee}{\Sigma{D}}
$$
satisfying $\rho{\comp}\sigma \homotopic 1_X$.
\end{thm}
\begin{proof}
Let $X$ be a finite complex and a co-H-space.
Then by \cite{Ganea:co-h-space}, $X$ is dominated by $\Sigma\Omega{X}$, where $\Sigma\Omega{X} \homotopic \Sigma\pi \vee \Sigma(\bigvee_{\tau \in \pi} \Omega_{\tau}X)$ and $\pi = \pi_1(X)$ is a free group.
As $X$ is a finite complex, the rank of $\pi$ is finite and $B = B\pi$ is a finite bunch of circles and $X$ is dominated by $B \vee \Sigma(\bigvee_{\tau \in \pi} \Omega_{\tau}X)$.
The image of $X$ can be taken as a finite subcomplex of $B \vee \Sigma(\bigvee_{\tau \in \pi} \Omega_{\tau}X)$ and so there is a finite subcomplex $D$ in $\bigvee_{\tau \in \pi} \Omega_{\tau}X$ such that $B{\vee}{\Sigma}D$ dominates $X$.

Conversely, let $X$ be dominated by $B{\vee}\Sigma{W}$ where $B$ is a finite bunch of circles indexed by a finite set $\Lambda$ and $W$ is a connected finite complex.
Then $B{\vee}\Sigma{W} = \Sigma(\Lambda{\vee}W)$, and $\cat{X} \leq \cat(\Sigma(\Lambda{\vee}W)) = 1$.
Thus $X$ is a co-H-space.
Also $X$ is dominated by the finite complex $B{\vee}\Sigma{W}$ whose fundamental group is free of finite rank.
The finiteness obstruction for $X$ lies in the Whitehead group $Wh(\pi)=K_0(\integer\pi_1(X))/\pm1$ (\cite{Wall:finiteness} and \cite{Varadarajan:wall-obstruction}) which is zero (see \cite{Cohn:free-ideal-ring} and \cite{Bass:projective}).
Thus $X$ has a homotopy type of a finite complex.
This completes the proof.
\end{proof}

Let $P = \sigma{\comp}\rho : B{\vee}{\Sigma{D}} \to B{\vee}{\Sigma{D}}$ be the homotopy idempotent given by Theorem \ref{thm:dominate}.
So $P$ restricted to $B$ can be chosen as the inclusion $\incl_B : B \subset B{\vee}{\Sigma{D}}$ and $P$ restricted to $\Sigma{D}$ lifts to $P_0 : {\Sigma{D}} \to \cover{B{\vee}{\Sigma{D}}}$ where $\cover{B{\vee}{\Sigma{D}}} \homotopic \bigvee_{\tau\in\pi} \tau\cdot{\Sigma{D}}$ as ${\Sigma{D}}$ is simply connected.
As ${\Sigma{D}}$ is a finite complex, $P_0({\Sigma{D}})$ is included in a finite subcomplex $\bigvee^{t}_{i=1} \tau_{i}\cdot{\Sigma{D}}$.
So the restriction of $P$ to $\Sigma{D}$ equals the composition
\begin{equation*}
{\Sigma{D}} \overset{P_0}\longrightarrow \begin{textstyle}\bigvee^{t}_{i=1}\end{textstyle} \tau_i\cdot{\Sigma{D}} \hookrightarrow \cover{B{\vee}{\Sigma{D}}} \overset{p}\longrightarrow B{\vee}{\Sigma{D}}.
\end{equation*}
Therefore we have the commutative diagram
\begin{equation}\label{cd:complete}
\begin{xy}
\xymatrix{
{B{\vee}{\Sigma{D}}}
	\ar[r]^{P}
	\ar[d]^{1_B{\vee}P_0}
&{B{\vee}{\Sigma{D}}}
\\
{B{\vee}\bigvee^{t}_{i=1} \tau_i\cdot{\Sigma{D}}}
	\ar@{^{(}->}[r]
&{B{\vee}\bigvee_{\tau\in\pi} \tau\cdot{\Sigma{D}}}
	\ar[u]^{\{\incl_B,p\}}
}\end{xy}
\end{equation}
which plays a crucial role in the next section.

We will need also maps $\sigma' : C \to \Sigma{D}$ and $\rho' : \Sigma{D} \to C$ defined by the commutative diagram below in which the columns are cofibrations:
\begin{equation}\label{diag:sigma_rho}
\begin{xy}
\xymatrix{
{B}
	\ar@{=}[r]
	\ar[d]
&{B}
	\ar@{=}[r]
	\ar[d]
&{B}
	\ar[d]
\\
{X}
	\ar[r]^{\sigma\quad}
	\ar[d]
&{B{\vee}\Sigma{D}}
	\ar[r]^{\quad\rho}
	\ar[d]
&{X}
	\ar[d]
\\
{C}
	\ar[r]^{\sigma'}
&{\Sigma{D}}
	\ar[r]^{\rho'}
&{C}
}
\end{xy}
\end{equation}
chosen so that $\rho'{\comp}\sigma' \homotopic 1_C$, as $\rho{\comp}\sigma \homotopic 1_X$.
Thus the self map $P'=\sigma'{\comp}\rho'$ of $\Sigma{D}$ is also a homotopy idempotent.
\section{Almost $p$-complete co-H-objects}
\par
Using the universality of almost $p$-completion, we have the natural equivalences between homotopy sets.
\begin{align*}&
[\apcomplete{(B{\vee}\Sigma{D})},\apcomplete{(B{\vee}\Sigma{D})}]_B
= [B{\vee}{\Sigma{D}},\apcomplete{(B{\vee}\Sigma{D})}]_B 
= [{\Sigma{D}},\apcomplete{(B{\vee}\Sigma{D})}] 
\\&\quad
= [{\Sigma{D}},\cover{\,\apcomplete{(B{\vee}\Sigma{D})}\,}] 
= [{\Sigma{D}},\pcomplete{(\cover{B{\vee}{\Sigma{D}}})}] 
%\\&\quad
%= [\pcomplete{\Sigma{D}},\pcomplete{\left(\begin{textstyle}\bigvee_{\tau{\in}\pi}\end{textstyle}\tau{\cdot}{\Sigma{D}}\right)}]
= [\pcomplete{\Sigma{D}},\begin{textstyle}\complete{\bigvee}_{\tau{\in}\pi}\end{textstyle}\tau{\cdot}\pcomplete{\Sigma{D}}]
\end{align*}
where $\complete{\bigvee}$ denotes the completed wedge sum.
% of $p$-complete spaces and $\sum$ implies `products of finitely many non-trivial elements'.
Projecting to its factors $\tau{\cdot}\pcomplete{\Sigma{D}}$, we have a map 
\begin{equation*}
\beta : [\apcomplete{(B{\vee}\Sigma{D})},\apcomplete{(B{\vee}\Sigma{D})}]_B
= [\pcomplete{\Sigma{D}},\begin{textstyle}{\complete{\bigvee}_{\tau{\in}\pi}}\end{textstyle}\tau{\cdot}\pcomplete{\Sigma{D}}] \rightarrow \begin{textstyle}\prod_{\tau{\in}\pi}\end{textstyle}\tau{\cdot}[\pcomplete{\Sigma{D}},\pcomplete{\Sigma{D}}]
\end{equation*}
to the product and the image of $\beta$ contains the sum $\begin{textstyle}\sum_{\tau{\in}\pi}\end{textstyle}\tau{\cdot}[\pcomplete{\Sigma{D}},\pcomplete{\Sigma{D}}]$. Indeed $\beta$ is a continuous homomorphism of topological groups, where the group structure is inherited from the co-H-space ${\Sigma {D}}$.  

We give an alternative description of the closed subgroup which is the
image of $\beta$.
Let $\{g_{\tau}\}_{\tau\in\pi}$ denote an element of the product $\begin{displaystyle}\prod_{\tau{\in}\pi}\end{displaystyle}\tau{\cdot}[\pcomplete{\Sigma{D}},\pcomplete{\Sigma{D}}]=\begin{displaystyle}\prod_{\tau{\in}\pi}\end{displaystyle}\tau{\cdot}[\Sigma{D},\pcomplete{\Sigma{D}}]$.
\begin{prop}\label{prop:beta}
The element $\{g_{\tau}\}_{\tau\in\pi}$ lies in the image of $\beta$ if and only if $\{\chi{\comp}g_{\tau}\}_{\tau\in\pi} \in \begin{textstyle}\sum_{\tau{\in}\pi}\end{textstyle}\tau{\cdot}[\Sigma{D},K]$ for any map $\chi : \pcomplete{\Sigma{D}} \to K$ and any space $K$ all of whose homotopy groups are finite p-groups.
\end{prop}
\begin{proof}
Let $f : \pcomplete{\Sigma {D}} \to \begin{textstyle}{\complete{\bigvee}_{\tau{\in}\pi}}\end{textstyle}\tau{\cdot}\pcomplete{\Sigma{D}}$ and $\beta(f) = \{f_{\tau}\}_{\tau\in\pi} \in \begin{displaystyle}\prod_{\tau{\in}\pi}\end{displaystyle}\tau{\cdot}[\Sigma{D},\pcomplete{\Sigma{D}}]$.
Since $(\begin{textstyle}{\complete{\bigvee}_{\tau{\in}\pi}}\end{textstyle}\tau{\cdot}\chi){\comp}f$ lies in $[\pcomplete{\Sigma {D}},\begin{textstyle}{\complete{\bigvee}_{\tau{\in}\pi}}\end{textstyle}\tau{\cdot}K] = [\Sigma {D},\begin{textstyle}{{\bigvee}_{\tau{\in}\pi}}\end{textstyle}\tau{\cdot}K]$, the map $\{\chi{\comp}f_{\tau}\}_{\tau\in\pi}$ lies in $\begin{textstyle}\sum_{\tau{\in}\pi}\end{textstyle}\tau{\cdot}[\Sigma{D},K]$ as required.
The converse statement holds by naturality and fundamental properties of $p$-completion.
\end{proof}
\par
\begin{lem}\label{lem:beta}
$\beta(\apcomplete{P}) \in \begin{textstyle}\sum_{\tau{\in}\pi}\end{textstyle}\tau{\cdot}[\pcomplete{\Sigma{D}},\pcomplete{\Sigma{D}}]$.
\end{lem}
\begin{proof}
The lemma follows from (\ref{cd:complete}).
\end{proof}
\par
We now recall results from \cite{Hubbuck:p-atomic} and \cite{Xu:thesis}.
We define a homomorphism of near-algebras by mapping homotopy classes
of self-maps of $\apcomplete{(B{\vee}\Sigma{D})}$ over $B$ to the induced endomorphism of $\tilde{H}_{\ast}(\cover{\,\apcomplete{(B{\vee}\Sigma{D})}\,};\field_p)
%\cong
%\tilde{H}_{\ast}(\pcomplete{\cover{(B{\vee}\Sigma{D})}};\field_p)
%\cong
%\field_p\pi{\otimes}\tilde{H}_{\ast}(\pcomplete{\Sigma{D}};\field_p)
\cong \tilde{H}_{\ast}(\pcomplete{\Sigma{D}};\field_p\pi)$ (see \cite{ISS:hom_co-h}) the $\field_p$-homology
groups of the universal cover
\begin{align*}&
\alpha : [\apcomplete{(B{\vee}\Sigma{D})},\apcomplete{(B{\vee}\Sigma{D})}]_B \to \End_{\field_p\pi}\{\tilde{H}_{\ast}(\pcomplete{\Sigma{D}};\field_p\pi)\}.
\end{align*}
When $B$ is a point, the same definition gives a homomorphism 
\begin{align*}&
\alpha_0 : [\pcomplete{\Sigma{D}},\pcomplete{\Sigma{D}}] \to \End_{\field_p}\{\tilde{H}_{\ast}(\pcomplete{\Sigma{D}};\field_p)\}.
\end{align*}
The homomorphisms $\alpha$ and $\alpha_0$ fit into a commutative diagram
\begin{equation*}
\begin{xy}
\xymatrix{
{[\apcomplete{(B{\vee}\Sigma{D})},\apcomplete{(B{\vee}\Sigma{D})}]_B}
\ar[rr]^{\alpha}
\ar[d]^{\beta}
&{}
&{\End_{\field_p\pi}\{\tilde{H}_{\ast}(\pcomplete{\Sigma{D}};\field_p\pi)\}}
\ar@{=}[dd]
\\
{\begin{displaystyle}\prod_{\tau{\in}\pi}\end{displaystyle}\tau{\cdot}[\pcomplete{\Sigma{D}},\pcomplete{\Sigma{D}}] }
\\
{\begin{displaystyle}\sum_{\tau{\in}\pi}\end{displaystyle}\tau{\cdot}[\pcomplete{\Sigma{D}},\pcomplete{\Sigma{D}}] }
\ar@{^{(}->}[u]
\ar[rr]^{\sum_{\tau{\in}\pi}\tau{\cdot}\alpha_0\qquad}
&{}
&{\begin{displaystyle}\sum_{\tau{\in}\pi}\end{displaystyle}\tau{\cdot}\End_{\field_p}\{\tilde{H}_{\ast}(\pcomplete{\Sigma{D}};\field_p)\}. }
}\end{xy}
\end{equation*}
Then $\alpha([\apcomplete{P}]) \in \End_{\field_p\pi}\{\tilde{H}_{\ast}(\pcomplete{\Sigma{D}};\field_p\pi)\}$, by Lemma \ref{lem:beta}.

The topological radical $N$ in the compact, Hausdorff space $[\pcomplete{\Sigma{D}},\pcomplete{\Sigma{D}}]$ is defined by
\begin{equation*}
%N = \{n \in [\pcomplete{\Sigma{D}},\pcomplete{\Sigma{D}}] \,\vert\, \mbox{$h{\comp}n$ is topologically nilpotent for all $h \in [\pcomplete{\Sigma{D}},\pcomplete{\Sigma{D}}]$}\}.
N = \{n \in [\pcomplete{\Sigma{D}},\pcomplete{\Sigma{D}}] \,\vert\, \mbox{$h{\comp}n$ is topologically nilpotent for all $h$}\}.
\end{equation*}
The radical $R$ in the finite ring $\End_{\field_p}\{\tilde{H}_{\ast}(\pcomplete{\Sigma{D}};\field_p)\}$ is defined by 
\begin{equation*}
R = \{r \in \End_{\field_p}\{\tilde{H}_{\ast}(\pcomplete{\Sigma{D}};\field_p)\} \,\vert\, \mbox{For any $u$, there is $v$ such that $v(1+ur) = 1$}\}.
\end{equation*}
Then (see section 3 in \cite{Hubbuck:p-atomic}), $\alpha_0$ induces a homomorphism of rings
\begin{equation*}
\alpha'_{0} : [\pcomplete{\Sigma{D}},\pcomplete{\Sigma{D}}]/N \rightarrow \left\{\End_{\field_p}\{\tilde{H}_{\ast}(\pcomplete{\Sigma{D}};\field_p)\}\right\}/R
\end{equation*}
which is a monomorphism onto its image, which can be identified with $\overset{k}{\underset{i=1}{\oplus}}M(n_i,\field_{q_i})$ for some $k$, where $\field_{q}$ is a finite field of characteristic $p$ with $q$ elements.
\begin{lem}\label{lem:p-atomic}
There is an isomorphism of rings induced by $\alpha$.
\begin{equation*}
\alpha' :  \begin{displaystyle}\sum_{\tau{\in}\pi}\end{displaystyle}\tau{\cdot}([\pcomplete{\Sigma{D}},\pcomplete{\Sigma{D}}]/N) \rightarrow \overset{k}{\underset{i=1}{\oplus}}M(n_i,\field_{q_i}\pi).
\end{equation*}
\end{lem}
\begin{proof}
We identify ${\field_p\pi}{\otimes}\End_{\field_p}\{\tilde{H}_{\ast}(\pcomplete{\Sigma{D}};\field_p)\}$ with $\End_{\field_p\pi}\{\tilde{H}_{\ast}(\pcomplete{\Sigma{D}};\field_p\pi)\}$ so that the image ${\field_p\pi}{\otimes}(\overset{k}{\underset{i=1}{\oplus}}M(n_i,\field_{q_i}))$ of $\alpha'=\sum_{\tau{\in}\pi}\tau{\cdot}\alpha_0$ becomes $\overset{k}{\underset{i=1}{\oplus}}M(n_i,\field_{q_i}\pi)$.
\end{proof}
\section{The proof of Theorem \ref{thm:main}}
\par
Lemmas \ref{lem:beta} and \ref{lem:p-atomic} imply that $\apcomplete{P} \in [\apcomplete{(B{\vee}\Sigma{D})},\apcomplete{(B{\vee}\Sigma{D})}]$ is mapped in homology to a direct sum of idempotents $\overset{k}{\underset{i=1}{\oplus}}P_i \in \overset{k}{\underset{i=1}{\oplus}}M(n_i,\field_{q_i}\pi)$.
We appeal to work of Bass \cite{Bass:k-theory}; each $P_i$ defines an $\field_{q_i}\pi$-homomorphism of $(\field_{q_i}\pi)^{n_i}$ and so there exists an $\field_{q_i}\pi$-isomorphism of $(\field_{q_i}\pi)^{n_i}$, $A_i$ say, such that
\begin{equation*}
\newcommand{\bg}{\fontfamily{cmr}\fontsize{21}{10pt}\selectfont}
\newcommand{\bigzero}{\smash{\hbox{\bg 0}}}
A_iP_iA^{-1}_i = %
	\begin{bmatrix}
	1 &          &   &   &          &   \\
	  & \ddots   &   &   & \bigzero &   \\
	  &          & 1 &   &          &   \\
	  &          &   & 0 &          &   \\
	  & \bigzero &   &   & \ddots   &   \\
	  &          &   &   &          & 0
	\end{bmatrix}%
		\in M(n_i,\field_{q_i}\pi).
\end{equation*}
This matrix also lies in $M(n_i,\field_{q_i})$.
Let a ring homomorphism $\epsilon : \field_q\pi \to \field_q$ be defined by $\epsilon(\tau) = 1$ and so $\epsilon(P_i)$ represents $\apcomplete{P'}$.
We choose $\phi : \apcomplete{(B{\vee}\Sigma{D})} \to \apcomplete{(B{\vee}\Sigma{D})}$ lying in $\beta^{-1}(\Sigma_{\tau\in\pi}\tau\cdot[\pcomplete{\Sigma{D}},\pcomplete{\Sigma{D}}])$ and a corresponding $\phi' : \pcomplete{\Sigma{D}} \to \pcomplete{\Sigma{D}}$ representing the invertible matrices $\overset{k}{\underset{i=1}{\oplus}}A_i$ in $M(n_i,\field_{q_i}\pi)$ and $\overset{k}{\underset{i=1}{\oplus}}\epsilon(A_i)$ in $M(n_i,\field_{q_i})$ respectively.
So $\phi$ and $\phi'$ are homotopy equivalences as they induce
isomorphisms of homology groups of universal covers by Lemma \ref{lem:p-atomic}.
Referring back to (\ref{diag:sigma_rho}), we have a commutative diagram
\begin{equation*}
\begin{xy}
\xymatrix{
{\apcomplete{X}}
	\ar[r]^{\apcomplete{\sigma}~\quad}
	\ar[d]
&{\apcomplete{(B{\vee}\Sigma{D})}}
	\ar[r]^{\phi}
	\ar[d]
&{\apcomplete{(B{\vee}\Sigma{D})}}
	\ar[r]^{~\quad\apcomplete{\rho}}
	\ar[d]
&{\apcomplete{X}}
	\ar[d]
\\
{\pcomplete{C}}
	\ar[r]^{\pcomplete{\sigma'}}
&{\pcomplete{\Sigma{D}}}
	\ar[r]^{\phi'}
&{\pcomplete{\Sigma{D}}}
	\ar[r]^{\pcomplete{\rho'}}
&{\pcomplete{C}.}
}
\end{xy}
\end{equation*}
\par
The self map $\phi{\comp}\apcomplete{P}{\comp}\phi^{-1}$ of $\apcomplete{(B{\vee}{\Sigma}D)}$ is a homotopy idempotent represented by the matrix $\overset{k}{\underset{i=1}{\oplus}}A_iP_iA^{-1}_i = \overset{k}{\underset{i=1}{\oplus}}\epsilon(A_i)\epsilon(P_i)\epsilon(A_i)^{-1}$ which also represents $\phi'{\comp}\pcomplete{P'}{\comp}{\phi'}^{-1}$.
Therefore this matrix also represents $\apcomplete{s{\comp}r}$ where $s=1_B{\vee}(\phi'{\comp}\pcomplete{\sigma'}) : B \vee \pcomplete{C} \to B \vee \pcomplete{{\Sigma}D}$ and $r=1_B{\vee}(\pcomplete{\rho'}{\comp}{\phi'}^{-1}) : B \vee \pcomplete{{\Sigma}D} \to B \vee \pcomplete{C}$, and so $r{\comp}s \homotopic 1_{B{\vee}\pcomplete{C}}$.
We deduce %that $\beta(\phi)$, $\beta(\apcomplete{P})$ and $\beta(\phi^{-1})$ are in $\begin{textstyle}\sum_{\tau\in\pi}\end{textstyle}\tau{\cdot}[\pcomplete{{\Sigma}D},\pcomplete{{\Sigma}D}]$ and satisfy
\begin{equation*}
\beta(\phi{\comp}\apcomplete{P}{\comp}\phi^{-1})
%= \beta(\phi{\comp}\apcomplete{P}{\comp}\phi^{-1}) 
\homotopic \beta(\apcomplete{(s{\comp}r)}) \mod \begin{textstyle}\sum_{\tau\in\pi}\end{textstyle}\tau{\cdot}N\quad\text{in}\quad\begin{textstyle}\sum_{\tau\in\pi}\end{textstyle}\tau{\cdot}[\pcomplete{{\Sigma}D},\pcomplete{{\Sigma}D}].
%\quad \text{in}\quad \begin{textstyle}\sum_{\tau\in\pi}\end{textstyle}\tau{\cdot}[\pcomplete{{\Sigma}D},\pcomplete{{\Sigma}D}].
\end{equation*}
\par
Let $f = \apcomplete{\rho}{\comp}\phi^{-1}{\comp}\apcomplete{s} : \apcomplete{(B{\vee}C)} \to \apcomplete{X}$ and $g = \apcomplete{r}{\comp}\phi{\comp}\apcomplete{\sigma} : \apcomplete{X} \to \apcomplete{(B{\vee}C)}$.
Then %$g{\comp}f : \apcomplete{(B{\vee}C)} \to \apcomplete{(B{\vee}C)}$ satisfies
\begin{align*}&
g{\comp}f 
	= (\apcomplete{r}{\comp}\phi{\comp}\apcomplete{\sigma}){\comp}(\apcomplete{\rho}{\comp}\phi^{-1}{\comp}\apcomplete{s})
%\\&
	= \apcomplete{r}{\comp}(\phi{\comp}\apcomplete{P}{\comp}\phi^{-1}){\comp}\apcomplete{s}
\intertext{whose image by $\beta$ is in $\begin{textstyle}\sum_{\tau\in\pi}\end{textstyle}\tau{\cdot}[\pcomplete{{\Sigma}D},\pcomplete{{\Sigma}D}]$ and is homotopic modulo $\begin{textstyle}\sum_{\tau\in\pi}\end{textstyle}\tau{\cdot}N$ to that of}
	%&\homotopic
&\apcomplete{r}{\comp}\apcomplete{(s{\comp}r)}{\comp}\apcomplete{s}
	= \apcomplete{(r{\comp}s)}{\comp}\apcomplete{(r{\comp}s)}
	%\quad\mod\begin{textstyle}\sum_{\tau\in\pi}\end{textstyle}\tau{\cdot}N
%\\&
	\homotopic \apcomplete{(1_{B{\vee}\pcomplete{C}})}{\comp}\apcomplete{(1_{B{\vee}\pcomplete{C}})}
	= 1_{\apcomplete{(B{\vee}C)}}.
\end{align*}
Thus the self map $g{\comp}f$ of $\apcomplete{(B{\vee}C)}$ over $B$ induces an isomorphisms of homology groups of the universal cover by Lemma \ref{lem:p-atomic}.
%$\cover{g{\comp}f} : \cover{\apcomplete{(B{\vee}C)}} \to \cover{\apcomplete{(B{\vee}C)}}$ is a $\field_p\pi$-homology equivalence of simply-connected $p$-complete spaces, and hence is a homotopy equivalence.
Therefore
\begin{enumerate}
\item[(4.1)]\quad
$g{\comp}f : \apcomplete{(B{\vee}C)} \to \apcomplete{(B{\vee}C)}$ is a homotopy equivalence.
\end{enumerate}
\par
It is routine to check that
\begin{enumerate}
\item[(4.2)]\quad
$f$ and $g$ induce isomorphisms of the $\field_p$-homology groups of universal covers,
\item[(4.3)]\quad
$f$ and $g$ induce isomorphisms of fundamental groups.
\end{enumerate}
\par
We complete the proof by following \cite{BD:ass-co-act}.
Statements (4.1), (4.2) and (4.3) are similar to the conclusion of Lemma 1.6 of \cite{BD:ass-co-act}.
One then makes minor changes to the proof of Theorem 1.5 given there to deduce Theorem \ref{thm:main}.
%
%     BIBLIOGRAPHY
%
\ifx\undefined\bysame
\newcommand{\bysame}{\leavevmode\hbox to3em{\hrulefill}\,}
\fi

\end{document}